\documentclass[%
amsmath,amssymb,
showkeys
]{revtex4-1}

\usepackage{graphicx}
\usepackage{dcolumn}
\usepackage{bm}
\usepackage{tikz}
\usepackage{float}
\usepackage{tabularx}

\usepackage{amsmath}
\usepackage{amsfonts}
\usepackage{amssymb}
\usepackage{graphicx}
\usepackage{verbatim}

\usepackage{lineno}

\usepackage{amsthm}
\newtheorem{definition}{Definition}
\newtheorem{theorem}{Theorem}
\newtheorem{lemma}{Lemma}[theorem]

\usepackage{hyperref}
\usepackage{environ}
\usepackage{algorithm}
\usepackage{algorithmic}

\begin{document}

\title{Strong contraction mapping and topological non-convex optimization}

\author{Siwei Luo}

\email{siuluosiwei@gmail.com}
\affiliation{%
The University of Illinois at Chicago, 1200 W Harrison St, Chicago, IL 60607
}%

\begin{abstract}

The strong contraction mapping, a self-mapping that the range is always a subset of the domain, admits a unique fixed-point which can be pinned down by the iteration of the mapping. We introduce a topological non-convex optimization method as an application of strong contraction mapping to achieve global minimum convergence. The strength of the approach is its robustness to local minima and initial point position. 







\end{abstract}



\maketitle


\newpage

\section{\label{sec:level1}Introduction}

Calculus of variation plays a critical role in the modern calculation. Mostly, the function of interest cannot be obtained directly but makes a functional attain its extremum. This leads people to construct different functionals according to different questions, such as the least action principle, Fermat's principle, maximum likelihood estimation, finite element analysis, machine learning and so forth. These methods provide us the routine to transfer the original problem to an optimization problem, which makes optimization methods have been used almost everywhere.

However, the popular gradient-based optimization methods be applied to a non-convex function with many local minima is unprincipled. They face great challenges of finding the global minimum point of the function. Because the information from the derivative of a single point is not sufficient to decide the global geometrical property of the function. In the gradient-based methods, the domain of the searching area is divided into several subsets with regards to different local minima. And usually, iterating points will converge to one local minimum depends on where the initial point starts from. For a successful minimum point convergence, the initial point luckily happens to be sufficiently near the global minimum point.

It is the very time to think outside box and try to cultivate some method other than gradient-based methods. Before that, let us first think about why gradient-based optimization methods fail the task of global minimum convergence. Let (X,d) be a metric space and T:X $\rightarrow$ X is a self-mapping. For the inequality that,
\begin{equation}
d(T(x),T(y)) \leq qd(x,y), \forall x,y \in X.
\end{equation}
if $q \in [0,1),$ T is called contractive; if $q \in [0,1],$ T is called nonexpansive; if $q < \infty $, T is called Lipschitz continuous\cite{Husain,Latif}. The gradient-based methods are usually nonexpansive mapping. The fixed point may exist but is not unique for a general situation. For instance, if the objective function $f(x)$ has many local minima, for the gradient descent method mapping $T(x) = x-\gamma \nabla f(x)$, any local minimum is a fixed point accordingly. From the perspective of spectra of a bounded operator, for a nonexpansive mapping, any minimum of the objective function is an eigenvector of the eigenvalue equation $T(x) = \lambda x$ , in which $\lambda = 1$. In the optimization problem, nonexpansive mapping sometimes works but their weakness is obvious. 

It is worth to note that to summarize the condition for the existence of fixed-point of non-expansive mapping or Lipschitz continuous mapping is a great challenge. Topologically, to prove the existence of fixed-point of the contraction mapping is relatively easier than to prove the existence of fixed-point of non-expansive mapping or Lipschitz continuous mapping. For contraction mapping, the set of the range is shrinking. While, for non-expansive mapping or Lipschitz continuous mapping, the change of range set can involve complicated transformation and rotation\cite{Ahues,Rudin}. There is no straightforward way to find the fixed-point of mapping. The fixed-point(or fixed-points) of Lipschitz continuous mapping, even if it exists, can cunningly hide inside an inflating volume. Because both the existence and uniqueness of solution are important in optimization problem so that the contractive mapping is more favored than the nonexpansive mapping.

The well-known Banach fixed-point theorem, as the first fixed-point theorem on contraction mapping, plays an important role in solving linear or nonlinear system. But for optimization problems, the condition of contraction mapping $T:X\rightarrow X$ that $d(T(x),T(y)) \leq qd(x,y)$, which usually requires convexity of the function, is hard to be applied to non-convex optimization problem. In this study, we are trying to extend the Banach fixed-point theorem to an applicable method for optimization problems, which is called strong contraction mapping. 

Strong contraction mapping is a surjection self-mapping that always maps to the subset of its domain. We will prove that strong contraction mapping admits a unique fixed-point, explain how to build an optimization method as an application of strong contraction mapping and illustrate why its fixed-point is the global minimum point of the objective function. 

\section{\label{sec:level1}Strong contraction mapping and the fixed-point}

Recall the definition of diameter $D(X)$ of a metric space $X$.
\begin{definition}
Let $(X,d)$ be a metric space. And the metric measurement $D(X)$ refers to the maximum distance between two points in the vector space X \cite{Fred} :
\begin{equation}
D(X) := sup \{ d(x,y), \forall x ,y \in X \}
\end{equation}\end{definition}

\begin{definition}
Let $(X,d)$ be a complete metric space. Then a mapping $T:X \rightarrow X$  is called weak contraction mapping on $X$ if the range of mapping T is always a subset of its domain during the iteration, namely, $\mathcal{R}(T) = X_{i+1} \subset X_i$ and there exists a $q\in [0,1)$ such that $D(X_{i+1}) \leq qD(X_i)$.
\end{definition} 


This contraction mapping is called strong because the requirement $D(X_{i+1}) \leq q(X_{i})$ is looser than $d(T(x),T(y)) \leq qd(x,y)$ what in Banach fixed-point theorem, which doesn't require the distance between two points getting smaller and smaller but the diameter of the range of the mapping getting smaller and eventually shrinking to a point. Thereafter, the inequality $d(T(x),T(y)) \leq qd(x,y)$ is included by the inequality $D(X_{i+1}) \leq q(X_{i})$ as a special case.

\begin{theorem}
Let $(X,d)$ be a non-empty complete metric space with strong contraction mapping $T:X \rightarrow X$. Then T admits a unique fixed-point $x^*$ in $X$ such that $x^* = T(x^*)$. 
\label{theorem_1}
\end{theorem}

To prove the Theorem.\ref{theorem_1}, one can follow the same logic of Banach fixed-point theorem proof\cite{Banach} but substitute the inequality $d(T(x),T(y)) \leq qd(x,y)$ with $D(X_{i+1}) \leq q(X_{i})$.  Let $x_0 \in X$ be arbitrary initial point and define a sequence $\{x_i\}$ as: $x_i = T(x_{i-1})$. 

\begin{lemma}
$\{x_i\}$ is a Cauchy sequence in $(X,d)$ and hence converges to a limit $x^*$ in $X$. \\
\end{lemma}
Proof.Let $m,n\in \mathbb{N}$ such that $m>n$.
\begin{equation*}
\begin{split}
	d(x_m,x_n) & \leq D(X_n) \\
    &\leq q^n D(X_0) \\
\end{split}
\end{equation*}
Let $\epsilon > 0$ be arbitrary, since $q \in [0,1)$, we can find a large $N\in \mathbb{N}$ such that
\begin{equation*}
	q^N \leq \frac{\epsilon }{D(X_0)}.
\end{equation*}
Hence, by choosing $m,n$ large enough: 
\begin{equation*}
	d(x_m,x_n) \leq q^nD(X_0) \leq \frac{\epsilon}{D(X_0)}D(X_0) = \epsilon.
\end{equation*}
Thus, $\{x_i\}$ is Cauchy sequence. $\Box$

As long as $D(X_0)$ is bounded, the convergence is guaranteed, which is independent of the choice of $x_0$.

\begin{lemma}
$x^*  :=  \lim x_{n}$ is a fixed-point of $T$ in $X$.\\
\end{lemma}
Proof. 
\begin{equation*}
\begin{split}
x^* = \lim_{i\rightarrow \infty}x_i = \lim_{n\rightarrow \infty}T(x_{i-1})\\
x^* = \lim_{i\rightarrow \infty}x_i = T(\lim_{n\rightarrow \infty}x_{i-1})
\end{split}	
\end{equation*}
Thus,$x^*=T(x^*)$. $\Box$

\begin{definition}
$X^* := \lim X_{i}$  is a fixed-set of T in X.
\end{definition}

\begin{lemma}
$x^*$ is the only fixed-point of $T$ in $ (X,d)$, the only element of $X^*$ and the diameter $D(X^*) = 0$ .\\
\label{lemma3}
\end{lemma}
Proof. Suppose there exists another fixed-point y that $T(y)=y$, then choose the subspace $X_i$ that both the $x^*$ and $y$ are the only elements in $X_i$. By definition, $X_{i+1} = \mathcal{R}(T(X_i))$ so that, both the $x^*$ and $y$ are elements in $X_{i+1}$, namely,

\begin{equation*}
  	0 \leq d(x^*,y) \leq D(X_{i+1}) \leq qD(X_i) = qd(x^*,y)  \\
\end{equation*}
\begin{equation*}
d(x^*,y) = 0
\end{equation*}
\begin{equation*}
 D(X^*) = D(\lim_{i\rightarrow \infty}X_{i}) = \lim_{i\rightarrow \infty}D(X_i) \leq \lim_{i\rightarrow \infty} q^n D(X_0) = 0 
\end{equation*} 

Since $(X,d)$ is a non-empty complete metric space and the diameter  $D(X^*) = 0$, $X^*$ has a single element $x^*$.  

Thus $x^* = y$. $\Box$

In this section, we have proven the existence and uniqueness of fixed-point of strong contraction mapping. Compared with contraction map, the strong contraction map can address a much wider range of problems as the requirement $D(T(x_i)) \leq D(x_i)$ is looser than $d(T(x),T(y)) \leq d(x,y)$. Different from $d(T(x),T(y)) \leq d(x,y)$, the inequality $D(T(x_i)) \leq D(x_i)$ doesn't require $x_i$ in sequence $\{x_n\}$ must move closer to each other for every step but confine the range of mapping to be smaller and smaller. Thereafter, the sequence $\{x_n\}$ is Cauchy sequence and converge to the fixed-point. It is worth to mention if some mapping $T$ is not a strong contraction mapping. Still, one can think about whether $T^n$ or $T^{-n}$ is a strong contraction mapping to find a way around. If that is the case, then the iteration of $T^n$ or $T^{-n}$ yields a fixed-point.



\section{\label{sec:level1}Optimization algorithm implementation}

After the discussion of strong contraction mapping, let us think about how to construct an optimization algorithm occupied with the property of strong contraction mapping. 

The objective function $f$ is a mapping defined on X, $f:X \rightarrow R$.

\begin{definition}
And an affine hyperplane is a subset $H$ of $X \times R$  of the form
\begin{equation*}
H = \{x \in X \times R; h(x) = L\}
\end{equation*}where $h$ is a linear functional that does not vanish identically and $L\in \mathcal{R}$ is a given constant\cite{Brezis}. 
\end{definition}

To overcome the dilemma that the optimization method may be saturated by local minima, intuitively, one can utilize a hyperplane that parallels to the domain to cut the objective function $f$ such that the intersection between the hyperplane and the function are contours. They reflect the global geometrical property of the function. The difficulty is how to iterate the hyperplane to move downwards and decide the position of the global minimum point. 

\begin{definition}
Contours is a subset C of $X \times R$ of the form 
\begin{equation*}
C = \{x\in X \times R;h(x) = L \in R, f(x) = L \in R \}
\end{equation*}\end{definition}

Observing that contours will divide the objective function $f$ into two parts. One is higher than the height of contours and the other is lower than the height of contours. Now, our task is to map to a point lower than the height of contours. As a numerical method, instead of getting all points of contours or symbolic expression of contours, we want to get n roots on the contours using root finding algorithm, that is 
\begin{equation}
    r_{j}^i = f^{-1}(L_i) = f^{-1}(f(x_{i})), \hspace{2pt} \forall j \in \{1,2,...,n\}
\end{equation}
where, $i$ indicates $i$th iterate and $j$ indicate $j$th root. 

First, provide one arbitrary initial point $x_{0}$ to the function and calculate the height of contours $f(x_{0}) = L_0$ at the height of the initial point; 

If the objective function $f$ is non-convex, as a consequence, to map to a point lower than the height of contours cannot be achieved by simply averaging the roots on contours. However, the objective function $f$ mostly can be decomposed into many locally convex subsets. We can map to a point lower than the height of contours by averaging the roots belong to the same locally convex subsets and iterate hyperplane downwards. 

A function $g:\mathcal{R}^n \rightarrow \mathcal{R}$ is convex, if for every $x,y \in \mathcal{R}^n$ and $\lambda \in [0,1]$, we have the inequality \cite{Bertsekas}
\begin{equation}
  g(\lambda x + (1 - \lambda) y ) \leq \lambda g(x) + (1-\lambda)g(y)
  \label{convex}
\end{equation}
The inequality $d(x_m,x_n) \leq q^n D(X_0)$ indicates the rate of convergence. The smaller q is, the higher rate of convergence is achieved. To achieve high rate of convergence, it is important to extend the Inequality.\ref{convex} to the situation with regards to many points $r_1,r_2,...,r_n$ , 
\begin{equation*}
  g(\lambda_{1}r_{1}+\lambda_{2}r_{2}+ ... +\lambda_{n}r_{n}) \leq \lambda_{1}g(r_1) + (1-\lambda_{1})g(\frac{\lambda_{2}r_{2}+ ... +\lambda_{n}r_{n}}{1-\lambda_{1}})
\end{equation*}
\begin{equation*}
  g(\lambda_{1}r_{1}+\lambda_{2}r_{2}+ ... +\lambda_{n}r_{n}) \leq \lambda_{1}g(r_1) + \lambda_{2}g(r_2) + (1-\lambda_{1}-\lambda_{2})g(\frac{\lambda_{3}r_{3}+ ... +\lambda_{n}r_{n}}{1-\lambda_{1}-\lambda_{2}})
\end{equation*}By induction,
\begin{equation}
  g(\lambda_{1}r_{1}+\lambda_{2}r_{2}+ ... +\lambda_{n}r_{n}) \leq \lambda_{1}g(r_1) + \lambda_{2}g(r_2) + ... + \lambda_{n}g(r_n)
\end{equation}
which is Jensen's Inequality \cite{Chandler},
\begin{equation}
  g(\sum_{j=1}^{n} \lambda_{j} r_{j}) \leq \sum_{j=1}^{n} \lambda_{j} g(r_{j})
  \label{Jensen}
\end{equation}
where, \begin{equation*}
    \sum_{i=1}^{n} \lambda_{i} = 1  
\end{equation*}
When we apply the Jensen's Inequality.\ref{Jensen} to the objective function $f$. Let $\lambda_{j} = \frac{1}{n}$,for $ j= 1, 2 , ... , n $,  the Jensen's Inequality.\ref{Jensen} turns to be a strict inequality, that is
\begin{equation}
  f(\lambda_{1}r_{1}+\lambda_{2}r_{2}+ ... +\lambda_{n}r_{n}) <   \lambda_{1}f(r_1) + \lambda_{2}f(r_2) + ... + \lambda_{n}f(r_n) = f(x_i) = L_i 
\end{equation}
Let $x_{i+1}  = \frac{1}{n} \sum_{j=1}^{n}r_{j}^{i}$  , that is
\begin{equation*}
 f(x_{i+1}) = L_{i+1} < f(x_i) = L_i, \hspace{10pt} \forall i = 1,2,...,n
\end{equation*}
Therefore, it is important to check whether two roots belong to the same locally convex subset and traverse all roots, decompose them into these locally convex subsets and get averages accordingly. Based on the Inequality.\ref{convex}, one practical way to achieve that is to pair every two roots on contours and scan function's value along the segment between the two roots and check whether there exists a point higher than contours' height. Traverse all the roots and apply this examination on them, then we can decompose the roots with regards to different locally convex subsets. To check whether two roots belong to the same convex subset, $N$ number of random points along the segment between two roots are checked whether higher than the contour's level or not\cite{Schachter,Tseng}. If the Inequality.\ref{convex} always holds for $N$ times of test, then we believe the two roots locate at the same locally convex subset and store them together. 

After the set of roots being decomposed into several locally convex subsets, the average of roots in the same locally convex subset is always lower than the contours' height due to Jensen's Inequality.\ref{Jensen}. 

\begin{theorem}
Provided there is a unique global minimum point $x_{min}$ of the objective function $f$, then the fixed-point $x^*$ of the strong contraction mapping $T$ is the global minimum point of the function.
\end{theorem}

Since the iterating point $x_{i+1}$ is always lower than $x_i$ ,
\begin{equation*}
\begin{split}
    0 \leq f(x_{i+1}) - f(x_{min}) < f(x_{i}) - f(x_{min}) \\
\end{split}
\end{equation*}Hence, there exists a $ \xi  \in (0,1) $ such that, 

\begin{equation*}
 0 \leq f(x_{i+1}) - f(x_{min}) < \xi^{i}(f(x_0) - f(x_{min})) 
\end{equation*}

As $i$ goes to infinity, then 
\begin{equation*}
\begin{split}
 0 \leq \lim_{i \rightarrow \infty}f(x_{i+1}) - f(x_{min}) &< \lim_{i \rightarrow \infty} \xi^{i}(f(x_0) - f(x_{min})) \\
  \lim_{i \rightarrow \infty} f(x_{i+1}) &- f(x_{min}) = 0 \\ 
    \lim_{i \rightarrow \infty} f(x_{i+1}) &= f(x_{min})\\
  f(\lim_{i \rightarrow \infty} x_{i+1}) &= f(x_{min}) \\
  f(x^*) &= f(x_{min})
\end{split}
\end{equation*}
Because the fixed-point $x^*$ is at the same height as the global minimum point $x_{min}$ and the global minimum point is unique. Thus, the fixed-point $x^*$ coincides with the global minimum point. The iteration $x_{i+1}=T{x_{i}}$ yields the fixed-point $x^*$that solves the equation $T(x^*)=x^*$.\cite{Taylor}

Rewrite the mapping that $x_{i+1} = T(x_i)$,which is averaging of roots $r_1,r_2,...r_n$ that locate at the same locally convex subset,
\begin{equation}
    x_{i+1} = T(x_i) = \frac{1}{n} \sum_{j=1}^{n}r_{j} = \frac{1}{n} \sum_{j=1}^{n}f^{-1}(f(x_{i}))
\end{equation}

Then, like a Russian doll, the dynamical system can be explicitly written as an expansion of iterates as
\begin{equation}
\begin{split}
x_{i+1} &= T(x_i) = T^{i+1}(x_0)\\
    &= \frac{1}{n}\sum_{j=1}^{n}r_{j} = \frac{1}{n}\sum_{j=1}^{n}f^{-1}(f(x_{i}))\\
    &= \frac{1}{n}\sum_{j=1}^{n}f^{-1}(f(\frac{1}{q} \sum_{j=1}^{q}f^{-1}(f(x_{i-1}))))\\
    &\vdots \\
    &= \frac{1}{n}\sum_{j=1}^{n}f^{-1}(f(\frac{1}{q} \sum_{j=1}^{q}f^{-1}(f( \dots \frac{1}{k} \sum_{j=1}^{k}f^{-1}(f(x_{0}))  \dots))))
\end{split}
\end{equation}
After the set of roots are decomposed into several convex subsets, the averages of roots with regards to each subset are calculated and the lowest one is returned as an update point from each iterate. Thereafter, the remaining calculation is to repeat the iterate over and over until convergence and return the converged point as the global minimum. The decomposition of roots set provided a divide-and-conquer method to transfer the original problem to a number of subproblems and solve them recursively. \cite{Cormen} 

In summary, the main procedure of topological non-convex optimization method is decomposed into the following steps: 1. Given the initial guess point $x_{0}$ for the objective function and calculate the contour level $L_0$;
2. Solve the equation $f(x_i)=L_i$ and get $n$ number of roots. Decompose the set of roots into several convex subsets, return the lowest average of roots as an update point from each iterate;
3. Repeat the above iterate until convergence. 

\begin{algorithm}[H] 
   \caption{topological non-convex optimization}
   \label{alg:weak_contraction_mapping}
\begin{algorithmic}
   \STATE {\bfseries input:}  $x_{0}$; 
   \STATE
   \STATE {\bfseries initialize:} set tolerance $\epsilon$; 
   \STATE
   \STATE {\bfseries repeat:}
   \STATE \hspace{10pt} calculate $L_i = f(x_{i})$;
   \STATE
   \STATE \hspace{10pt} $vector{<}root{>}$  $convex[Num\_Of\_Convex] $;
   \STATE
    \STATE {\bfseries \hspace{10pt} for}$(j=1;j<m;j{++})$ 
   \STATE \hspace{40pt} $r_{j}^i = f^{-1}(L_i) = f^{-1}(f(x_{i}))$; \hspace{10pt} $\Diamond$ solved by root finding algorithm
   \STATE
   \STATE {\bfseries \hspace{10pt} Traverse and pair two roots $r_{p}^i$ and $r_{q}^{i}$}
   \STATE
   \STATE \hspace{40pt} bool flag = true;
   \STATE
   \STATE {\bfseries \hspace{40pt} for}$(k=1;k<N;k{++})$\{
   \STATE \hspace{70pt}  $\lambda$ = random(0,1); 
   \STATE {\bfseries \hspace{70pt} if}$(f(\lambda r_{p}^i+(1-\lambda)r_{q}^i)>=L_i)$ {flag = false,\bfseries break}; \hspace{10pt} $\Diamond$ check whether they belong to the same locally convex subset
   \STATE \hspace{40pt} \}
   \STATE
   \STATE {\bfseries \hspace{40pt} if} (flag)\{
   \STATE \hspace{70pt} $convex[\zeta].push\_back(r_{p}^i)$; 
   \STATE \hspace{70pt} $convex[\zeta].push\_back(r_{q}^i)$;  \hspace{10pt} $\Diamond$ store roots into a same locally convex subset
   \STATE {\bfseries \hspace{40pt}} \}
   \STATE
   \STATE {\bfseries \hspace{10pt} for}$(\zeta=0;\zeta<convex.size();\zeta{++})$  
   \STATE \hspace{40pt} $\alpha_{\zeta}= \frac{1}{n} \sum_{j=1}^{n} r_{j}^i$; 
   \STATE
   \STATE \hspace{10pt} $x_{i+1} = argmin (f(\alpha_{\zeta}))$; \hspace{10pt} $\Diamond$ choose the lowest average to achieve higher rate of convergence
   \STATE
   \STATE {\bfseries until:} {$x_{i+1} - x_{i} < \epsilon$}
   \STATE
   \STATE {\bfseries return:} $x_{i+1}$
\end{algorithmic}
\end{algorithm}

According to the previous discussion, the mechanism of the proposed optimization method has posted some restrictions to the objective function. In summary, it requests the objective function to be continuous, occupied with a unique global minimum point, can be decomposed into many locally convex components and its domain is non-empty complete. And the good news is the most of real-world non-convex optimization problems can satisfy these prerequisites.

\section{\label{HSA}Experiments on Sphere, McCormick and Ackley functions}

\begin{figure}[H]
\includegraphics[width=1\textwidth]{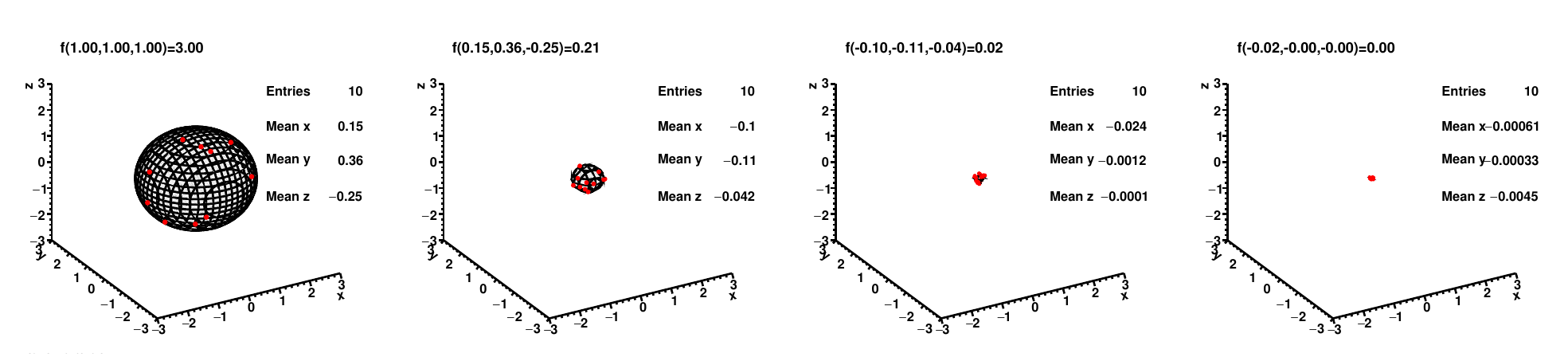}\centering
\caption{The red point markers are the roots and spherical surface is the contour is 3D space for each iteration. The contour in 3D space is a ball for Sphere function. The intermediate steps illustrate how the contour shrink to a point during the procedure.}
\label{Sphere_function}
\end{figure}

First of all, the optimization algorithm has been tested on a convex function the Sphere function $f(\boldmath{x})=\sum{x^{2}_{i}}$. The minimum is (0,0,0), where $f(0,0,0)=0$. The iterations of roots and contours is shown in FIG.\ref{Sphere_function} and the update of searching point is shown in TABLE.\ref{Sphere_table}.

\begin{table}[H]
\begin{center}
\begin{tabular*}{0.5\textwidth}{ c  c  c}
\textbf{iterate} & \space{       }\space{     } \textbf{ updating point} & \space{ }\space{     }\space{}\space{} \textbf{height of contour} \\
\hline\hline
0 & (1.00,1.00,1.00) & 3.0000\\ \centering
1 & (0.15,0.36,-0.25) & 0.2146\\ \centering
2 & (-0.1,-0.11,-0.042) & 0.0237 \\
3 & (-0.024,-0.0012,-0.0001) & 0.0004 \\
4 & (-0.0061,-0.00033,-0.0045) & \\
\hline
\end{tabular*}
\end{center}
\caption{When the optimization method is tested on Sphere function, the average of roots and the level of contour for each iteration is shown above.}\centering
\label{Sphere_table}
\end{table}

\begin{figure}[H]
\includegraphics[width=0.5\textwidth]{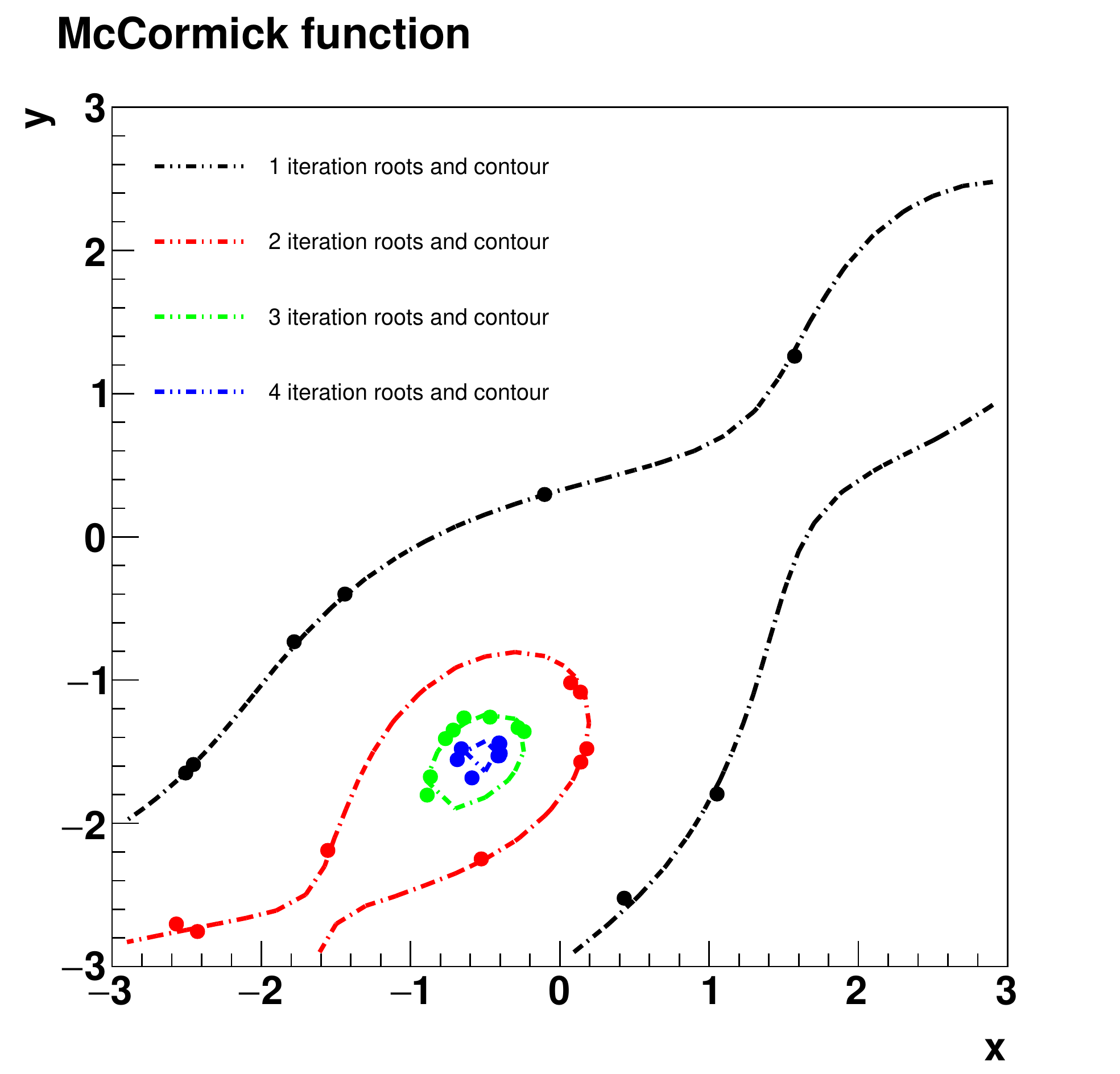}\centering
\caption{The point markers are the roots and the dash lines are the contours during iterates. The first 4 iteration results are drawn to illustrate the procedure of optimization test on the McCormick function.}
\label{McCormick_function}
\end{figure}

Then, we test the optimization algorithm on the McCormick function. And the first 4 iterates of roots and contour is shown in FIG.\ref{McCormick_function} and the detailed iteration of the searching point from the numerical calculation is shown in TABLE.\ref{McCormick_table}. The test result indicate the average of roots can move towards the global minimum point (-0.54719,-1.54719), where $f(-0.54719,-1.54719) = -1.9133$. 

\begin{table}[H]
\begin{center}
\begin{tabular*}{0.5\textwidth}{ c  c  c}
\textbf{iterate} & \space{       }\space{     } \textbf{updating point} & \space{ }\space{     }\space{}\space{} \textbf{height of contour} \\
\hline\hline
0 & \space{} (2.00000,2.0000) & 2.2431975047\\ \centering
1 & \space{} (-0.6409,-0.8826) & -1.1857067055\\ \centering
2 & \space{} (-0.8073,-1.8803) & -1.7770492114\\
3 & \space{} (-0.5962,-1.4248) & -1.8814760998\\
4 & \space{} (-0.4785,-1.5162) & -1.9074191216\\
5 & \space{} (-0.5640,-1.5686) & -1.9125755974\\
6 & \space{} (-0.5561,-1.5467) & -1.9131043354\\
7 & \space{} (-0.5474,-1.5465) & -1.9132219834\\
8 & \space{} (-0.5473,-1.5472) & \\
\hline

\end{tabular*}
\end{center}
\caption{When the optimization method is tested on McCormick function, the average of roots and the level of contour for each iteration is shown above.}
\label{McCormick_table}
\end{table}

\begin{figure}[H]
\includegraphics[width=0.5\textwidth]{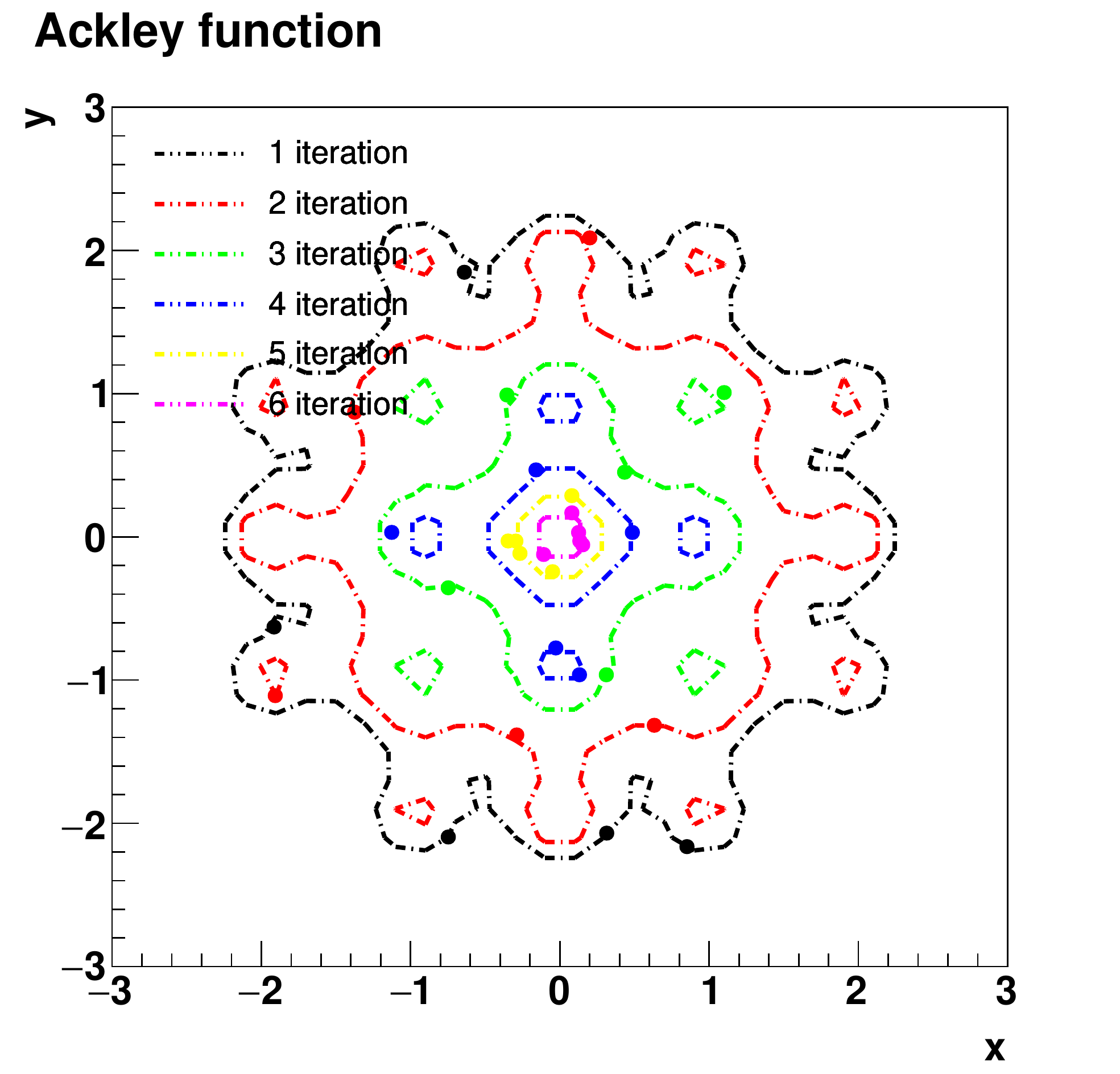}\centering
\caption{The point markers are the roots and the dash lines are contours during iterates. The first 6 iteration results are drawn to illustrate the procedure of optimization test on the Ackley function. The intermediate steps show how contours ignore the existence of local minima and safely approach to the global minimum point.}
\label{Ackley_function}
\end{figure}

Ackley function would be a nightmare for most gradient-based methods. Nevertheless, the algorithm has been tested on Ackley function where the global minimum locates at (0,0) that $f(0,0)=0$. And the first 6 iterates of roots and contours are shown in FIG.\ref{Ackley_function} and the minimum point (-0.00000034,0.00000003) return by the algorithm is shown in TABLE.\ref{Ackley_table}. The test result shows that the optimization algorithm is robust to local minima and able to achieve the global minimum convergence. The quest to locate the global minimum pays off handsomely. 

\begin{table}[H]
\begin{center}
\begin{tabular*}{0.5\textwidth}{ c  c  c}
\textbf{iterate} & \space{       }\space{     } \textbf{updating point} & \space{ }\space{     }\space{}\space{} \textbf{height of contour} \\
\hline\hline
0 & \space{} (2.00000000,2.00000000) & 6.59359908\\ \centering
1 & \space{} (-0.78076083,-1.34128187) & 5.82036224\\ \centering
2 & \space{} (-0.35105371,-0.62030933) & 4.11933422\\
3 & \space{} (-0.20087095,0.38105138) & 3.09359564\\
4 & \space{} (0.06032320,-0.88101860) & 2.17077104\\
\vdots & \vdots & \vdots\\
15 & \space{} (0.00000404,-0.00000130) & 0.00001199\\
16 & \space{} (-0.00000194,-0.00000079) &0.00000591\\
17 & \space{} (-0.00000034,0.00000003) &\\
\hline
\end{tabular*}
\end{center}
\caption{When the optimization method is tested on Ackley function, the average of roots and the level of contour for each iteration is shown above.}
\label{Ackley_table}
\end{table}

The observation from these experiments is that the size of contour become smaller and smaller during the iterative process and eventually converge to a point, which is the global minimum point of the function.

\section{\label{sec:level1}Conclusion}

We introduced the definition of the strong contraction mapping and the existence and uniqueness of its fixed-point in this paper. As an extension of Banach fixed-point theorem, the iteration of strong contraction mapping is a Cauchy sequence and yields the unique fixed-point, which perfectly fit with the task of optimization. The global minimum convergence regardless of local minima and initial point position is a very significant strength for the optimization algorithm. And we illustrated how to implement an optimization method occupied with strong contraction mapping property. This topological optimization method finds a way around that even if the objective function is non-convex, still, we can decompose it into many convex components and take advantage of convexity of a locally convex component to pin down the global minimum point. This optimization method has been tested on Sphere, McCormick and Ackley functions and successfully achieved the global minimum convergence as expected. These experiments demonstrate the contours' shrinking and the iterating point's approaching the global minimum point. We look forward to extending our study to the higher dimensional situation and believe that the optimization method works for that in principle.


\begin{thebibliography}{9}

\bibitem{Husain}
Husain, T., and Abdul Latif. "Fixed points of multivalued nonexpansive maps." International Journal of Mathematics and Mathematical Sciences 14.3 (1991): 421-430.

\bibitem{Latif}
Latif, Abdul, and Ian Tweddle. "On multivalued f-nonexpansive maps." Demonstratio Mathematica 32.3 (1999): 565-574.

\bibitem{Ahues}
Ahues, Mario, Alain Largillier, and Balmohan Limaye. Spectral computations for bounded operators. Chapman and Hall/CRC, 2001.

\bibitem{Rudin}
Rudin, Walter. "Functional analysis. 1991." Internat. Ser. Pure Appl. Math (1991).

\bibitem{Brezis}
Brezis, Haim. Functional analysis, Sobolev spaces and partial differential equations. Springer Science Business Media, 2010.

\bibitem{Fred}
Croom, Fred H. Principles of topology. Courier Dover Publications, 2016.

\bibitem{Kiwiel}
Kiwiel, K. C. (2001). Convergence and efficiency of subgradient methods for quasiconvex minimization. Mathematical programming, 90(1), 1-25.

\bibitem{Banach}
Banach, S. (1922). Sur les opérations dans les ensembles abstraits et leur application aux équations intégrales. Fund math, 3(1), 133-181.

\bibitem{Branciari}
Branciari, A. (2002). A fixed point theorem for mappings satisfying a general contractive condition of integral type. International Journal of Mathematics and Mathematical Sciences, 29(9), .

\bibitem{Taylor}
Taylor, A. E., \& Lay, D. C. (1958).  (Vol. 2). New York: Wiley.

\bibitem{Schachter}
Schachter, Bruce. "Decomposition of polygons into convex sets." IEEE Transactions on Computers 11 (1978): 1078-1082.

\bibitem{Tseng}
Tseng, Paul. "Applications of a splitting algorithm to decomposition in convex programming and variational inequalities." SIAM Journal on Control and Optimization 29.1 (1991): 119-138.

\bibitem{Arveson}
Arveson, W. (2006). A short course on spectral theory (Vol. 209). Springer Science  Business Media.

\bibitem{Blanchard}
Blanchard, P.; Devaney, R. L.; Hall, G. R. (2006). Differential Equations. London: Thompson. pp. 96–111. ISBN 0-495-01265-3.

\bibitem{Bertsekas}
Bertsekas, Dimitri (2003). Convex Analysis and Optimization. Athena Scientific.

\bibitem{Chandler}
Chandler, David. "Introduction to modern statistical mechanics." Introduction to Modern Statistical Mechanics, by David Chandler, pp. 288. Foreword by David Chandler. Oxford University Press, Sep 1987. ISBN-10: 0195042778. ISBN-13: 9780195042771 (1987): 288.

\bibitem{Cormen}
Cormen, Thomas H., Charles E. Leiserson, Ronald L. Rivest, and Clifford Stein. Introduction to algorithms. MIT press, 2009.

\end{thebibliography}
\end{document}